\input amssym.def
\input epsf

\let \blskip = \baselineskip
\parskip=1.2ex plus .2ex minus .1ex

\tabskip 20pt
\tolerance = 1000
\pretolerance = 50
\newcount\itemnum
\itemnum = 0
\overfullrule = 0pt

\def\title#1{\bigskip\centerline{\bigbigbf#1}}
\def\author#1{\bigskip\centerline{\bf #1}\smallskip}
\def\address#1{\centerline{\it#1}}
\def\abstract#1{\vskip1truecm{\narrower\noindent{\bf Abstract.} #1\bigskip}}

\def\sp{\bigskip}
\def\nosp{\vskip -\the\blskip plus 1pt minus 1pt}

\def\br{\hfil\break} 
\def\ti{\br \hglue \the \parindent}

\def\skipit#1{}
\def\mdag{\raise 3pt\hbox{\dag}}
\def\mddag{\raise 3pt\hbox{\ddag}}
\def\mstar{\raise 3pt\hbox{$\ast$}}

\def\XP{\par\noindent\hang}
\def\LP{\par\noindent}
\def\BP[#1]{\par\item{[#1]}}
\def\SH#1{\sp\vskip\parskip\leftline{\bigbf #1}\nobreak}

\def\TH#1{\sp\XP{\bf THEOREM\ \shead#1}}
\def\LM#1{\sp\XP{\bf LEMMA\ \shead#1}}

\def\CO#1{\sp\XP{\bf COROLLARY\ \shead#1}}

\def\PF{\LP{\bf Proof:\ }}

\def\NX{\advance\itemnum by 1 \sp\LP {\bf \shead \the\itemnum.\ }}
\def\qed{\null\nobreak\hfill\hbox{${\vrule width 5pt height 6pt}$}\par\sp}
\def\qedno{\null\nobreak\hfill\hbox{${\vrule width 5pt height 6pt}$}}
\def\eqed{\eqno{\hbox{${\vrule width 5pt height 6pt}$}}}

\def\cart{\>\hbox{${\vcenter{\vbox{
    \hrule height 0.4pt\hbox{\vrule width 0.4pt height 4.5pt
    \kern4pt\vrule width 0.4pt}\hrule height 0.4pt}}}$}\>}
\def\bxmu{\>\hbox{${\vcenter{\vbox {
    \hrule height 0.4pt\hbox{\vrule width 0.4pt height 4pt
    \hskip -1.3pt\lower 1.8pt\hbox{$\times$}\negthinspace\vrule width 0.4pt}
    \hrule height 0.4pt}}}$}\>}

\def\lin#1{\hbox to #1true in{\hrulefill}}



\def\JGT{{\it J.\ Graph Theory}}

\def\GnC{{\it Graphs and Combin.}}

\def\CNum{{\it Congr.\ Numer.{}}}

\def\al{\alpha}			
    
	 	\def\DLT{\Delta}

\def\NN{{\Bbb N}}  
 \def\ZZ{{\Bbb Z}}




\def\({\left(}	\def\){\right)}


\def\CL#1{\left\lceil{#1}\right\rceil}   \def\CFR#1#2{\CL{{#1\over#2}}}

\def\st{\colon\,}   

\def\SET#1:#2{\{#1\colon\;#2\}}

		



\magnification=\magstep1
\vsize=9.0 true in
\hsize=6.5 true in
\headline={\hfil\ifnum\pageno=1\else\folio\fi\hfil}
\footline={\hfil\ifnum\pageno=1\folio\else\fi\hfil}

\parindent=20pt
\baselineskip=12pt
\parskip=.2ex  

\def\shead{ }

\font\bigbf = cmb10 scaled \magstep1

\font\bigbigbf = cmb10 scaled \magstep2


\def\gpic#1{#1
     \medskip\par\noindent{\centerline{\box\graph}} \medskip}

\title{COLORING OF TREES WITH MINIMUM SUM OF COLORS}
\author{Tao Jiang and Douglas B. West}
\address{j-tao@math.uiuc.edu and west@math.uiuc.edu}
\address{University of Illinois, Urbana, IL 61801-2975}
\vfootnote{}{\br
   Running head: COLORING OF TREES \br
}
\abstract{The {\it chromatic sum} $\Sigma(G)$ of a graph $G$ is the smallest
sum of colors among all proper colorings with natural numbers.  The {\it
strength} $s(G)$ of $G$ is the minimum number of colors needed to achieve the
chromatic sum.  We construct for each positive integer $k$ a tree $T_k$  with
strength $k$ that has maximum degree only $2k-2$.  The result is best possible.
}

\SH
{1. INTRODUCTION}

A {\it proper coloring} of the vertices of a graph $G$ is a function
$f\colon\;V(G)\to{\Bbb N}$ such that adjacent vertices receive different
labels (colors).  The {\it chromatic number} $\chi(G)$ is the minimum number of
colors in a proper coloring of $G$.  The {\it chromatic sum} $\Sigma(G)$ is a
variation introduced by Ewa Kubicka in her dissertation.  It is the
minimum of $\sum_{v\in V(G)} f(v)$ over proper colorings $f$ of $G$.
A {\it minimal coloring} of $G$ is a proper coloring of
$G$ such that $\sum_v f(v) = \Sigma(G)$.

One might think that a minimal coloring can be obtained by selecting a proper
coloring with the minimum number of colors and then giving the largest color
class color $1$, the next largest color $2$, and so on.  However, even among
trees, which have chromatic number $2$, more colors may be needed to obtain a
minimal coloring.  The {\it strength} $s(G)$ of a graph $G$ is the minimum
number of colors needed to obtain a minimal coloring.  Kubicka and Schwenk [4]
constructed for every positive integer $k\geq 2$ a tree $T_k$ with strength $k$.
Thus $s(G)$ may be arbitrarily large even when $\chi(G)=2$
(trivially $s(G) \geq \chi(G)$).

How large can $s(G)$ be in terms of other parameters?
When vertices are colored greedily in natural numbers with respect to
a vertex ordering $v_1,\ldots,v_n$, the number of colors used
is at most $1+\max_i d^*(v_i)$, where $d^*(v_i)$ counts the
neighbors of $v_i$ in $\{v_1,\ldots,v_{i-1}\}$.  Always this yields
$\chi(G)\le 1+\Delta(G)$.  The best upper bound on $\chi(G)$ that can be
obtained in this way is the
{\it Szekeres-Wilf number} $w(G)=1+\max_{H\subseteq G}\delta(H)$
(also confusingly called the ``coloring number'').
Interestingly, the average of these two well-known upper bounds for the
chromatic number is an upper bound for the strength $s(G)$.

\TH{(Hajiabolhassan, Mehrabadi, and Tusserkani [2])} Every
graph $G$ has strength at most $\CL{(w(G)+\Delta(G))/ 2}$.

\sp
We show that this bound is sharp, even for trees.
Every nontrivial tree $T$ has Szekeres-Wilf number $2$, and thus
$ s(T) \leq 1+{\lceil{\Delta(T)/ 2}\rceil}$.
In the Kubicka-Schwenk construction [4], the tree with strength $k$
has maximum degree about $k^2/2$.  To show that the bound above is sharp,
we construct for each $k\ge1$ a tree $T_k$ with strength $k$ and
maximum degree $2k-2$.  Given a proper coloring $f$ of a tree $T$, we use
$\Sigma f$ to denote $\sum_{v\in V(T)}f(v)$.

\SH
{2. THE CONSTRUCTION}

Linearly order the pairs of natural numbers so that $(h,l)<(i,j)$ if either
$h+l<i+j$ or $h+l=i+j$ and $l<j$.  With respect to this ordering, we
inductively construct for each pair $(i,j)\in{\Bbb N}\times{\Bbb N}$ a rooted
tree $T_i^j$ and a coloring $f_i^j$ of $T_i^j$.  In other words, we construct
trees in the order $T_1^1$, $T_2^1$, $T_1^2$, $T_3^1,\ldots$.  Our desired tree
with strength $k$ will be $T_k^1$.  Let $[n]=\{k\in\ZZ\st 1\le k\le n\}$.

\bigskip
\noindent
{\bf Construction.}
Let $T_1^1$ be a tree of order 1, and let $f_1^1$ assign color 1 to this
single vertex.  Consider $(i,j)\ne (1,1)$, and suppose that for each
$(h,l)<(i,j)$ we have constructed $T_h^l$ and $f_h^l$.  We construct $T_i^j$
and $f_i^j$ as follows.  Let $u$ be the root of $T_i^j$.  For each $k$ such that
$1\leq k \leq i+j-1$ and $k\neq i$, we take two copies of $T_k^m$, where
$m=\CL{(i+j-k)/2}$, and we let the roots of these $2(i+j-2)$ trees
be children of $u$.  The resulting tree is $T_i^j$ (see Fig.~1). 
Define the coloring $f_i^j$ of $T_i^j$ by assigning $i$ to the root $u$ and
using $f_k^m$ on each copy of $T_k^m$ rooted at a child of $u$.
\qed
\def\To{T_1^{\CFR{i+j-1}2}}
\def\Tk{T_k^{\CFR{i+j-k}2}}
\def\Tn{T_{i+j-1}^1}

\gpic{
\expandafter\ifx\csname graph\endcsname\relax \csname newbox\endcsname\graph\fi
\expandafter\ifx\csname graphtemp\endcsname\relax \csname newdimen\endcsname\graphtemp\fi
\setbox\graph=\vtop{\vskip 0pt\hbox{%
    \special{pn 8}%
    \special{ar 298 1469 298 689 0 6.28319}%
    \special{ar 987 1469 298 689 0 6.28319}%
    \special{ar 1906 1469 298 689 0 6.28319}%
    \special{ar 2594 1469 298 689 0 6.28319}%
    \special{ar 3513 1469 298 689 0 6.28319}%
    \special{ar 4202 1469 298 689 0 6.28319}%
    \graphtemp=.5ex\advance\graphtemp by 0.781in
    \rlap{\kern 0.298in\lower\graphtemp\hbox to 0pt{\hss $\bullet$\hss}}%
    \graphtemp=.5ex\advance\graphtemp by 0.781in
    \rlap{\kern 0.987in\lower\graphtemp\hbox to 0pt{\hss $\bullet$\hss}}%
    \graphtemp=.5ex\advance\graphtemp by 0.781in
    \rlap{\kern 1.906in\lower\graphtemp\hbox to 0pt{\hss $\bullet$\hss}}%
    \graphtemp=.5ex\advance\graphtemp by 0.781in
    \rlap{\kern 2.594in\lower\graphtemp\hbox to 0pt{\hss $\bullet$\hss}}%
    \graphtemp=.5ex\advance\graphtemp by 0.781in
    \rlap{\kern 3.513in\lower\graphtemp\hbox to 0pt{\hss $\bullet$\hss}}%
    \graphtemp=.5ex\advance\graphtemp by 0.781in
    \rlap{\kern 4.202in\lower\graphtemp\hbox to 0pt{\hss $\bullet$\hss}}%
    \graphtemp=.5ex\advance\graphtemp by 0.092in
    \rlap{\kern 2.250in\lower\graphtemp\hbox to 0pt{\hss $\bullet$\hss}}%
    \special{pa 298 781}%
    \special{pa 2250 92}%
    \special{pa 987 781}%
    \special{pa 2250 92}%
    \special{pa 1906 781}%
    \special{fp}%
    \special{pa 2594 781}%
    \special{pa 2250 92}%
    \special{pa 3513 781}%
    \special{pa 2250 92}%
    \special{pa 4202 781}%
    \special{fp}%
    \graphtemp=.5ex\advance\graphtemp by 0.000in
    \rlap{\kern 2.250in\lower\graphtemp\hbox to 0pt{\hss $u$\hss}}%
    \graphtemp=.5ex\advance\graphtemp by 0.781in
    \rlap{\kern 1.446in\lower\graphtemp\hbox to 0pt{\hss $\cdots$\hss}}%
    \graphtemp=.5ex\advance\graphtemp by 0.781in
    \rlap{\kern 3.054in\lower\graphtemp\hbox to 0pt{\hss $\cdots$\hss}}%
    \graphtemp=.5ex\advance\graphtemp by 1.469in
    \rlap{\kern 1.446in\lower\graphtemp\hbox to 0pt{\hss $\cdots$\hss}}%
    \graphtemp=.5ex\advance\graphtemp by 1.469in
    \rlap{\kern 3.054in\lower\graphtemp\hbox to 0pt{\hss $\cdots$\hss}}%
    \graphtemp=.5ex\advance\graphtemp by 1.469in
    \rlap{\kern 0.298in\lower\graphtemp\hbox to 0pt{\hss $\To$\hss}}%
    \graphtemp=.5ex\advance\graphtemp by 1.469in
    \rlap{\kern 0.987in\lower\graphtemp\hbox to 0pt{\hss $\To$\hss}}%
    \graphtemp=.5ex\advance\graphtemp by 1.469in
    \rlap{\kern 1.906in\lower\graphtemp\hbox to 0pt{\hss $\Tk$\hss}}%
    \graphtemp=.5ex\advance\graphtemp by 1.469in
    \rlap{\kern 2.594in\lower\graphtemp\hbox to 0pt{\hss $\Tk$\hss}}%
    \graphtemp=.5ex\advance\graphtemp by 1.469in
    \rlap{\kern 3.513in\lower\graphtemp\hbox to 0pt{\hss $\Tn$\hss}}%
    \graphtemp=.5ex\advance\graphtemp by 1.469in
    \rlap{\kern 4.202in\lower\graphtemp\hbox to 0pt{\hss $\Tn$\hss}}%
    \graphtemp=.5ex\advance\graphtemp by 2.296in
    \rlap{\kern 0.643in\lower\graphtemp\hbox to 0pt{\hss 2 copies\hss}}%
    \graphtemp=.5ex\advance\graphtemp by 2.296in
    \rlap{\kern 2.250in\lower\graphtemp\hbox to 0pt{\hss 2 copies\hss}}%
    \graphtemp=.5ex\advance\graphtemp by 2.296in
    \rlap{\kern 3.857in\lower\graphtemp\hbox to 0pt{\hss 2 copies\hss}}%
    \hbox{\vrule depth2.480in width0pt height 0pt}%
    \kern 4.500in
  }%
}%
}
\centerline{$1\leq k \leq i+j-1$ and $k \neq i$}
\bigskip
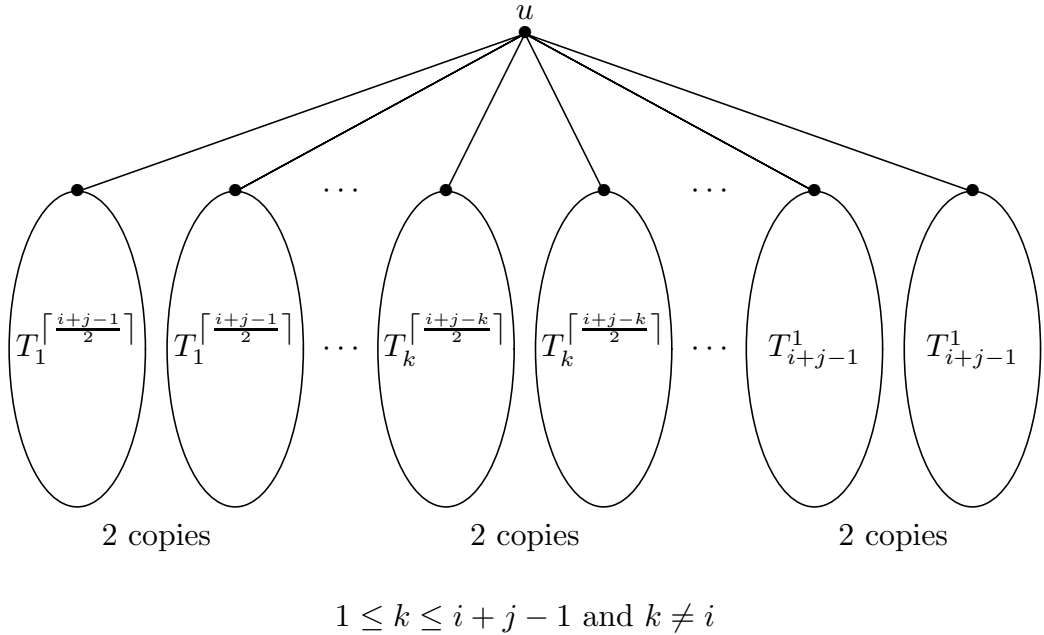
\centerline {Figure 1. The  construction of $T_i^j$ }

\LM{}
For $(i,j)\in \NN\times\NN$, the construction of $T_i^j$ is well-defined,
and $f_i^j$ is a proper coloring of $T_i^j$ with color $i$ at the root.
\PF
To show that $T_i^j$ is well-defined,
it suffices to show that when $(i,j)\ne(1,1)$, every tree used in the
construction of $T_i^j$ has been constructed previously.  We use trees
of the form $T_k^m$, where $k\in[i+j-1]-\{i\}$ and $m=\CL{(i+j-k)/2}$.
It suffices to show that $k+m\le i+j$ and that $m<j$ when $k+m=i+j$.

For the first statement, we have $k+m\le \CL{(i+j+k)/2}\le i+j$,
since $k\le i+j-1$.  Equality requires $k=i+j-1$, which occurs only when
$j\ge2$ and yields $m=1$.  Thus $m<j$ when $k+m=i+j$.  Since the trees whose
indices sum to $i+j$ are generated in the order
$T_{i+j-1}^1,\ldots,T_1^{i+j-1}$, the tree
$T_k^m$ exists when we need it.

Finally, $f_i^j$ uses color $i$ at the root of $T_i^j$, by construction.
Since the subtrees used as descendants of the root have the form
$T_k^m$ with $k\ne i$, by induction the coloring $f_i^j$ is proper.
\qedno

\SH
{3. THE PROOF}
The two-parameter construction enables us to prove a technically stronger
statement.  The additional properties of the construction facilitate the
inductive proof.  Recall that all colorings considered are labelings with
positive integers.

\TH{}
The construction of $T_i^j$ and $f_i^j$ has the following properties:

(1) If $f'$ is a coloring of $T_i^j$ different from $f_i^j$, then
$\Sigma f'>\Sigma f_i^j$.  Furthermore, if $f'$ assigns a color different from
$i$ to the root of $T_i^j$, then $\Sigma f' -\Sigma f_i^j \geq j$;

(2) If $j=1$, then $\Delta(T_i^j)=2i-2$, achieved by the root of $T_i^j$. 
If $j\geq 2$, then $\Delta(T_i^j)=2(i+j)-3$; 

(3) The highest color used in $f_i^j$ is $i+j-1$.
\PF
We use induction through the order in which the trees are constructed.
As the basis step, $T_1^1$ is just a single vertex, and $f_1^1$ gives it color
1; conditions (1)-(3) are all satisfied.

Now consider $(i,j)\ne(1,1)$.  For simplicity, we write $T$ for $T_i^j$ and $f$
for $f_i^j$.  To verify (1), let $f'$ be a coloring of $T$ different from $f$.
We consider two cases. 

{\bf Case 1. } $f'$ assigns $i$ to the root $u$ of $T$. 

In this case, $f'$ and $f$ differ on $T-u$.  Recall that $T-u$ is the
union of $2(i+j-2)$ previously-constructed trees.  The colorings $f'$
and $f$ differ on at least one of these trees.  By the induction hypothesis,
the total under $f'$ is at least the total under $f$ on each of these subtrees,
and it is larger on at least one.  Hence $\Sigma f'> \Sigma f$.

{\bf Case 2.} $f'$ assigns a color different from $i$ to the root $u$.

In this case, we need to show that $\Sigma f' - \Sigma f \geq j$. 
Again the induction hypothesis gives $f'$ as large a total as $f$ on
each component of $T-u$.  If $f'(u)\ge i+j$, then the difference on $u$
is large enough to yield $\Sigma f'-\Sigma f\ge j$.

Hence we may assume that $f'(u)=k$, where $1\leq k \leq i+j-1$ and $k\neq i$.
Since $f'$ is a proper coloring, it assigns a label other than $k$ to the roots
$v,v'$ of the two copies of $T_k^m$ in $T-u$, where $m=\CL{(i+j-k)/2}$.
Since $f$ uses $f_k^m$ on each copy of $T_k^m$, we have $f(v)=f(v')=k$.
Since $f'(v)$ and $f'(v')$ differ from $k$, the induction hypothesis implies
that on each copy of $T_k^m$ the total of $f'$ exceeds the total of $f$ by at
least $m$.  Since the total is at least as large on all other components,
we have
$$\Sigma f' - \Sigma f \geq k-i+2m = k-i+2\CL{i+j-k\over 2}\ge j.$$

Next we verify (2).  In the construction of $T=T_i^j$, we place $2(i+j-2)$
subtrees under the root $u$.  These have the form $T_k^m$ for 
$1\le k \le i-1$ and $i+1\le k\le i+j-1$, and always $m=\CL{(i+j-k)/2}$.
Note that $m=1$ only when $k=i+j-1$ or $k=i+j-2$.  The subtrees have
maximum degree $2k-2$ (when $m=1$) or $2(k+m)-3$ (when $m>1$).
Note that $2(k+m)-3>2k-2$ when $m\ge1$.  Thus
$$\Delta(T_k^m)\le 2(k+m)-3=2\left(k+\CL{i+j-k\over 2}\right)-3=
2\CL{i+j+k\over 2}-3.$$
Also, we always have $k+m=\CL{(i+j+k)/2}$ for the subtree $T_k^m$.

When $j=1$ we only have $k\le i-1$, and thus
$\Delta(T_k^m)\leq 2\CL{(i+1+k)/ 2}-3\leq 2i-3$.
Hence each vertex in $T-u$ has degree at most $(2i-3)+1=2i-2$ in $T$.
Since $d_T(u)=2i-2$, we have $\Delta(T)=2i-2$, achieved by the root.
  
When $j\geq 2$, the values of $k$ for the subtrees are $1\le k\le i-1$
and $i+1\le k\le i+j-1$.  By the induction hypothesis, the maximum degree of
$T_{i+j-1}^1$ is $2(i+j-1)-2=2(i+j)-4$ and is achieved by its root.  In $T$
this vertex has degree $2(i+j)-3$, which exceeds $d_T(u)$.
For $k\le i+j-2$, we have
$\Delta(T_k^m)\le 2\CL{(i+j+k)/ 2}-3\leq 2(i+j)-5$. 
Hence $\Delta(T)=2(i+j)-3$, achieved by the roots of the trees that
are isomorphic to $T_{i+j-1}^1$.

It remains to verify (3): the maximum color used in $f_i^j$ is $i+j-1$.
By the induction hypothesis and the construction,
the maximum color used by $f_k^m$ on each $T_k^m$ within $f_i^j$ is
$k+m-1=\CL{(i+j+k)/ 2}-1$.  Since the largest $k$ is $i+j-1$
when $j\ge2$ and is $i-1$ when $j=1$, this computation yields $i+j-1$ when
$j\ge2$ and $i-1$ when $j=1$ as the maximum color on $T-u$.
Since $f$ assigns $i$ to the root $u$, we obtain $i+j-1$ as the maximum color
on $T$ for both $j\ge2$ and $j=1$.
\qed

We have proved that $f_i^j$ is the unique minimal coloring of $T_i^j$ and
that it uses $i+j-1$ colors.  Hence $s(T_i^j)=i+j-1$. 
The maximum degree is $2i-2$ or $2(i+j)-3$, depending on whether $j=1$ or
$j\geq2$.  In particular, $T_i^1$ is a tree with strength $i$ and maximum degree
$2i-2$.

\CO{1.}
There exists for each positive integer $i$ a tree $T_i$ with
$s(T_i)=i$ and $\Delta(T_i)=2i-2$.
\qedno

\CO{2.}
For every real number $\alpha\in (0,1/2)$, there is a sequence of trees
$T_1',T_2',\ldots$ such that $\lim_{n\to \infty} {s(T_n')}/{\DLT(T_n')}=\al$.
\PF
Let $t=\lfloor({1\over \alpha}-2)i\rfloor+2$.  Consider the construction of
$T_i^1$.  Form $T_i'$ by adding $t$ additional
copies of the subtree $T_{i-1}^1$ under the root $u$ of $T_i^1$.
The strength of $T_i'$ is $i$, but $\Delta(T_i')=2i-2+t$.
As $i\rightarrow\infty$, we have 
$$ {s(T_i')\over \Delta(T_i')} = { i\over {2i+t-2}}={i\over {2i+ 
\lfloor({1\over \alpha}-2)i\rfloor}} \rightarrow \alpha.\eqed$$

\SH
{References}
\frenchspacing

\BP[1]
P. Erd\H os, E. Kubicka, and A.J. Schwenk,
Graphs that require many colors to achieve their chromatic sum,
\CNum\ 
71(1990), 17--28.

\BP[2]
H.~Hajiabolhassan, M.L. Mehrabadi, and R. Tusserkani,
Minimal coloring and strength of graphs.  {\it Proc.~ 28th Annual Iranian
Math. Conf., Part 1 (Tabriz, 1997)}, Tabriz Univ. Ser. 377 (Tabriz Univ.,
Tabriz, 1997), 353--357.

\BP[3]
E. Kubicka,
Constraints on the chromatic sequence for trees and graphs,
\CNum\ 76(1990), 219--230.

\BP[4]
E. Kubicka and A.J. Schwenk,
An introduction to chromatic sums,
{\it Proc. ACM Computer Science Conference, Louisville(Kentucky)}
1989, 39--45.

\BP[5]
C. Thomassen, P. Erd\H os, Y. Alavi, P.J. Malde, and A.J. Schwenk,
Tight bounds on the chromatic sum of a connected graph,
\JGT\ 13(1989), 353--357.

\BP[6]
Z. Tuza, Contraction and minimal $k$-colorability
\GnC\ 6(1990), 51--59.

\bye